\documentclass{elsarticle}
\usepackage{amssymb}
\usepackage{amsthm}
\usepackage{amsmath}
\usepackage{color}
\usepackage{enumerate}

\newtheorem{thm}{Theorem}
\newtheorem{lem}[thm]{Lemma}
\newtheorem{rem}[thm]{Remark}
\newtheorem{cor}[thm]{Corollary}
\newtheorem{defi}[thm]{Definition}
\newtheorem{prop}[thm]{Proposition}

\begin{document}

\title{$\zeta-$function and heat kernel formulae}

\author[syd]{Fedor Sukochev\corref{cor1}}
\ead{f.sukochev@unsw.edu.au}
\author[adl]{Dmitrii Zanin \fnref{arc}}
\ead{zani0005@csem.flinders.edu.au}

\cortext[cor1]{Corresponding Author}  
\fntext[arc]{Research supported by the
Australian Research Council}

\address[syd]{School of Mathematics and Statistics, University of New South Wales, Sydney, 2052, Australia.}
\address[adl]{School of Computer Science, Engineering and Mathematics, Flinders University, Bedford Park, 5042, Australia.}

\begin{abstract} We present a systematic study of asymptotic behavior of (generalised) $\zeta-$func\-tions and heat kernels used in noncommutative geometry and clarify their connections with Dixmier traces. We strengthen and complete a number of results from the recent literature and answer (in the affirmative) the question raised by M. Benameur and T. Fack \cite{BF}.

\end{abstract}


\maketitle

\section{Introduction}\label{introdu}

The interplay between Dixmier traces, $\zeta-$functions and heat kernel formulae is a cornerstone of noncommutative geometry \cite{Connes}. These formulae are widely used in physical applications. To define these objects, let us fix a Hilbert space $H$ and let $B(H)$ be the algebra of all bounded operators on $H$ with its standard trace ${\rm Tr}.$ Let $A$ and $B$ be positive operators from $B(H).$ Consider the following $[0,\infty]$-valued functions
\begin{equation}\label{zeta-def}
t\to\frac1t{\rm Tr}(A^{1+1/t}),\quad t\to\frac1t{\rm Tr}(A^{1+1/t}B))
\end{equation}
and, for fixed $0<q<\infty$
\begin{equation}\label{hk-q-def}
t\to\frac1t{\rm Tr}(\exp(-(tA)^{-q})),\quad t\to\frac1t{\rm Tr}(\exp(-(tA)^{-q})B).
\end{equation}
When these functions are finitely valued, they are frequently referred to as $\zeta-$functions and heat kernel functions associated with the operators $A$ and $B.$ When these functions are bounded, a particular interest is attached to their asymptotic behavior when $t\to\infty,$ which is usually measured with the help of some generalised limit $\gamma:L_{\infty}(0,\infty)\to\mathbb{R}$ yielding the following functionals
\begin{equation}\label{gamma-zeta-def}
\zeta_\gamma(A):=\gamma(\frac1t{\rm Tr}(A^{1+1/t})),\quad \zeta_{\gamma,B}(A):=\gamma(\frac1t{\rm Tr}(A^{1+1/t}B))
\end{equation}
and,
\begin{equation}\label{gamma-hk-q-def}
\varphi_\gamma(A):=\gamma(\frac1t{\rm Tr}(\exp(-(tA)^{-q}))),\quad \varphi_{\gamma,B}(A):=\gamma(\frac1t{\rm Tr}(\exp(-(tA)^{-q}))B).
\end{equation}
A natural class of operators for which the formulae \eqref{zeta-def} and \eqref{gamma-zeta-def} are well defined (respectively, \eqref{hk-q-def} and \eqref{gamma-hk-q-def}) is given by the set $\mathcal{M}_{1,\infty}$ (respectively, $\mathcal{L}_{1,\infty}$) of compact operators from $B(H).$ More precisely, denote by $\mu_{n}(T),$ $n\in\mathbb{N},$ the singular values of a compact operator $T$ (the singular values are the eigenvalues of the operator $|T|=(T^*T)^{1/2}$ arranged with multiplicity in decreasing order, (\cite[\S 1]{Simon}). Then
\begin{equation} \label{eq:dixmier_ideal}
\mathcal{M}_{1,\infty}:=\mathcal{M}_{1,\infty}(H)=\{T:\ \sup_{n\in\mathbb{N}}\frac1{\log(n+1)}\sum_{k=1}^n\mu_k(T)<\infty\}.
\end{equation}
defines a Banach ideal of compact operators. We set
$$\mathcal{L}_{1,\infty}:=\{T\in\mathcal{M}_{1,\infty}:\ \exists C>0\ \text{such that}\ \mu_n(A)\leq C/n,\ n\ge 1\}.$$
It is important to observe that the subset $\mathcal{L}_{1,\infty}$ is not dense in $\mathcal{M}_{1,\infty}$ (see e.g. \cite{KS}). It should also be pointed out that our notation here differs from that used in \cite{Connes}.

It follows from \cite[Theorem 4.5]{CRSS}  that the functions defined in \eqref{zeta-def} are bounded if and only if $A\in\mathcal{M}_{1,\infty}.$ It also follows from \cite{CRSS} and \cite{CGRS} that the functions defined in \eqref{hk-q-def} are bounded if and only if $A\in\mathcal{L}_{1,\infty}.$ In fact the last result is a strong motivation to consider the following modification of formulae \eqref{hk-q-def}. Let us consider a Cesaro operator on $L_{\infty}(0,\infty)$ given by
$$(Mx)(t)=\frac1{\log(t)}\int_1^tx(s)\frac{ds}{s},\quad t\in(0,\infty).$$
It follows from \cite{CRSS} and \cite{CGRS} that the functions
\begin{equation}\label{modified-hk-q-def}
M(t\to\frac1t{\rm Tr}(\exp(-(tA)^{-q}))),\quad M(t\to\frac1t{\rm Tr}(\exp(-(tA)^{-q}))B)
\end{equation}
are bounded if and only if $A\in\mathcal{M}_{1,\infty}.$ Therefore, for a given generalised limit $\omega,$ let us set
\begin{equation}\label{Connes states}
\omega':=\omega\circ M
\end{equation}
and instead of the functions given in \eqref{gamma-hk-q-def} consider the functions
\begin{equation}\label{modified-gamma-hk-q-def}
\xi_\omega(A):=\omega'(\frac1t{\rm Tr}(\exp(-(tA)^{-q}))),\quad \xi_{\omega,B}(A):=\omega'(\frac1t{\rm Tr}(\exp(-(tA)^{-q}))B).
\end{equation}

The class of dilation invariant states $\omega'$ as above was introduced by A. Connes (see \cite{Connes}) and it is natural to refer to this class as ''Connes states''. We prove in section \ref{heatk lin} that if $\omega$ in \eqref{Connes states} is dilation invariant, then $\xi_\omega$ is a linear functional on $\mathcal{M}_{1,\infty}.$ In fact, we also show in Proposition \ref{coincidence} that if $\omega$ in \eqref{Connes states} is such that $\xi_\omega$ is linear on $\mathcal{M}_{1,\infty},$ then
necessarily there exists a dilation invariant generalised limit $\omega_0$ such that $\xi_\omega=\xi_{\omega_0}.$

There is a deep reason to require that the functionals $\xi_\omega$ and $\zeta_\gamma$ be defined on $\mathcal{M}_{1,\infty}$ and be linear (and thus, by implication, to consider Connes states). Important formulae in noncommutative geometry \cite{Connes} and its semifinite counterpart \cite{CPS,CaS,BF,CRSS,CGRS} then connect these functionals with Dixmier traces on $\mathcal{M}_{1,\infty}.$ Recall that in \cite{Dixmier}, J. Dixmier constructed a non-normal semifinite trace  (a Dixmier trace) on $B(H)$ using the weight
\begin{equation} \label{eq:dtr0}
{\rm Tr}_\omega(T):=\omega\left(\left\{\frac{1}{\log(1+n)}\sum_{k=1}^n\mu_k(T)\right\}_{n=1}^\infty \right)\quad T>0,
\end{equation}
where $\omega$ is a dilation invariant state on $L_\infty(0,\infty).$

The interplay between positive functionals ${\rm Tr}_\omega,$ $\zeta_\gamma$ and $\xi_\omega$ on $\mathcal{M}_{1,\infty}$ makes
an important chapter in noncommutative geometry and has been treated (among many other papers) in \cite{Connes,CPS,CaS,BF,CRSS,S,CGRS, SSZ}. We now list a few most important known results concerning this interplay and explain our contribution to this topic.

In \cite{CPS}, the equality
\begin{equation}\label{cps zeta}
{\rm Tr}_{\omega}(AB)=(\omega\circ\log)(\frac1t\tau(A^{1+1/t}B))=\zeta_{\omega\circ\log,B}(A),\quad 0\leq A\in \mathcal{M}_{1,\infty}
\end{equation}
was established for every $B\in B(H)$ under very restrictive conditions on $\omega.$ These conditions are dilation invariance for both $\omega$ and $\omega\circ\log$ and $M-$invariance of $\omega.$ In \cite{CRSS}, for the special case $B=1,$ the assumption that $\omega$ is $M-$invariant has been removed. However, the case of an arbitrary $B$ appears to be inaccessible by the methods in that article. In Section \ref{zeta}, we prove the general result which implies, in particular, that the equality \eqref{cps zeta} holds without requiring $M-$invariance of $\omega.$

In \cite{CPS}, the equality
\begin{equation}\label{cps heat}
\omega(\frac1t\tau(\exp(-(tA)^{-q})B))=\Gamma(1+\frac1q)\tau_{\omega}(AB)
\end{equation}
was established under the same conditions on $\omega$ and $\omega\circ\log$ as above. In \cite{SSZ}, in the special case $B=1$ the equality \eqref{cps heat} was established under the assumption that $\omega$ is $M-$invariant. However, again the case of an arbitrary $B$ appears to be inaccessible by the methods in that article. Here, we are able to treat the case of a general operator $B.$

In \cite{BF} a more general approach to the heat kernel formulae is suggested. It consists of replacing the function $t\to\exp(t^{-q})$ with an arbitrary function $f$ from the Schwartz class. The following equality was proved in \cite{BF}
\begin{equation}\label{bf heat}
\omega(\frac1t\tau(f(tA)B))=\int_0^{\infty}f(\frac1s)ds\cdot\tau_{\omega}(AB)
\end{equation}
for $A\in\mathcal{L}_{1,\infty}$ and $M-$invariant $\omega.$

In \cite[p.51]{BF}, M. Benameur and T. Fack have asked whether the result above continues to stand without the $M-$invariance assumption on $\omega.$ In Theorem \ref{xihk red} below, we answer this question affirmatively for a much larger class of functions
than the Schwartz class and for any $A\in \mathcal{M}_{1,\infty}.$

Finally, it is important to emphasize the connection between our results with the theory of fully symmetric functionals. Recall that a linear positive functional $\varphi:\mathcal{M}_{1,\infty}\to\mathbb{C}$ is called {\it fully symmetric} if $\varphi(B)\leq\varphi(A)$ for every positive $A,B\in\mathcal{M}_{1,\infty}$ such that $B\prec\prec A.$ The latter symbol means that
$$\sum_{k=1}^n\mu_k(B)\leq \sum_{k=1}^n\mu_k(A),\quad \forall n\in\mathbb{N}.$$

It is obvious that every Dixmier trace ${\rm Tr}_\omega$ is a fully symmetric functional. However, the fact that every fully
symmetric functional coincides with a Dixmier trace is far from being trivial (see \cite{KSS} and Theorem \ref{fs is dixmier} below). It is therefore quite natural to ask whether a similar result holds for the sets of all linear positive functionals on $\mathcal{M}_{1,\infty}$ formed by the $\xi_\omega$ and $\zeta_\gamma$ respectively. To this end, we establish results somewhat similar to those of \cite{KSS}. Firstly, in Theorem \ref{xiomega linearity} we prove that if $\omega$ in \eqref{Connes states} is dilation invariant, then the functional $\xi_\omega$ extends to a fully symmetric functional on $\mathcal{M}_{1,\infty}.$ Secondly, in Theorem \ref{also main} we show that in fact every normalized fully symmetric functional on $\mathcal{M}_{1,\infty}$ coincides with some $\xi_\omega,$ where $\omega$ is dilation invariant. Thus, in view of \cite{KSS}, we can conclude that the set $\{{\rm Tr}_\omega:\ \omega\ \text{is a dilation invariant generalised limit}\}$ coincides with the set $\{\xi_\omega:\ \omega\ \text{is a dilation invariant generalised limit}\}$ (up to a norming constant). At the same time, a natural question, namely, whether the equality
$$\xi_\omega=\Gamma(1+\frac1q){\rm Tr}_{\omega}$$
holds for every dilation invariant generalised limit $\omega$ is answered in the negative in Theorem \ref{counterexample}.

Finally, we note that the question on the relationship between the sets $\{{\rm Tr}_\omega:\ \omega\ \text{is a dilation\ invariant\
generalised\ limit}\},$ $\{\zeta_\gamma:\ \gamma\ \text{is a generalised limit}\}$ and $\{\zeta_\omega:\ \omega\ \text{is\ a\ dilation\ invariant\ generalised\ limit}\}$ remains open.

\section{Definitions and notations}\label{definot}

The theory of singular traces on operator ideals rests on some classical analysis which we now review for completeness.

As usual, $L_{\infty}(0,\infty)$ is the set of all bounded Lebesgue measurable functions on the semi-axis equipped with the uniform norm $\|\cdot\|.$ Given a function $x\in L_{\infty}(0,\infty),$ one defines its decreasing rearrangement $\mu(x)=\mu(\cdot,x)$ by the formula (see e.g. \cite{KPS})
$$\mu(t,x)=\inf\{s\geq0:\ m(\{|x|>s\})\leq t\}.$$

Let $H$ be a Hilbert space and let $B(H)$ be the algebra of all bounded operators on $H$ equipped with the uniform norm $\|\cdot\|.$ Let $\mathcal{N}\subset B(H)$ be a semi-finite von Neumann algebra with a fixed faithful and normal semi-finite trace $\tau.$ For every $A\in\mathcal{N},$ the generalised singular value function $\mu(A)=\mu(\cdot,A)$ is defined by the formula (see e.g. \cite{FK})
$$\mu(t,A):=\inf\{\|Ap\|:\ \tau(1-p)\leq t\}.$$
If, in particular, $\mathcal{N}=B(H),$ then $\mu(A)$ is a step function and, therefore, can be identified with the sequence $\{\mu(n,A)\}_{n\geq0}$ of singular numbers of the operators $A$ (the singular values are the eigenvalues of the operator $|A|=(A^*A)^{1/2}$ arranged with multiplicity in decreasing order).

Equivalently, $\mu(A)$ can be defined in terms of the distribution function $d_A$ of $A.$ That is, setting
$$d_A(s):=\tau(e^{|A|}(s,\infty)),\ s\ge 0,$$
we obtain
$$\mu(t,A)=\inf\{s:\ d_A(s)\geq t\}, \ t>0.$$
Here, $e^{|A|}$ denotes the spectral measure of the operator $|A|.$

The following formula follows directly from the von Neumann definition of trace (see the definition at \cite[Definition 15.1.1]{vN})
\begin{equation}\label{f formula}
\tau(f(A))=-\int_0^{\infty}f(\lambda)dd_A(\lambda).
\end{equation}

Using the Jordan decomposition, every operator $A\in B(H)$ can be uniquely written as
$$A=(\Re(A)_+-\Re(A)_-)+i(\Im(A)_+-\Im(A)_-).$$
Here, $\Re(A):=\frac{1}{2}(A+A^*)$ (respectively, $\Im(A):=\frac{1}{2i}(A-A^*)$) for any operator $A\in B(H)$ and $B_+=Be^B(0,\infty)$ (respectively, $B_-=Be^B(-\infty,0)$) for any self-adjoint operator $B\in B(H).$ Recall that $\Re A,\Im A\in\mathcal{N}$ for every $A\in\mathcal{N}$ and $B_+,B_-\in\mathcal{N}$ for every self-adjoint $B\in\mathcal{N}.$

Let $\psi:\mathbb{R}_+\to\mathbb{R}_+$ be an increasing concave function such that $\psi(t)=O(t) $ as $t\to 0.$ The Marcinkiewicz function space $M_{\psi}$ (see e.g. \cite{KPS}) consists of all $x\in L_{\infty}(0,\infty)$ satisfying
$$\|x\|_{M_{\psi}}:=\sup_{t>0}\frac1{\psi(t)}\int_0^t\mu(s,x)ds<\infty.$$
The Marcinkiewicz operator space $\mathcal{M}_{\psi}:= \mathcal{M}_{\psi}(\mathcal{N},\tau)$ (see e.g. \cite{CaS,CRSS}) consists of all $A\in\mathcal{N}$ satisfying
$$\|A\|_{\mathcal{M}_{\psi}}:=\sup_{t>0}\frac1{\psi(t)}\int_0^t\mu(s,A)ds<\infty.$$
We are especially interested in Marcinkiewicz spaces $\mathcal{M}_{1,\infty}$ and $M_{1,\infty}$ that arise when $\psi(t)=\log(1+t),$ $t\ge 0.$ In the literature, the ideal $\mathcal{M}_{1,\infty}$ is sometimes referred to as the Dixmier ideal. We recommend the recent paper of A. Pietsch, \cite{Pie}, discussing the origin of $\mathcal{M}_{1,\infty}$ in mathematics.

For $s>0,$ dilation operators $\sigma_s:L_{\infty}\to L_{\infty}$ are defined by the formula $(\sigma_sx)(t)=x(t/s).$ Clearly, $\sigma_s:M_{1,\infty}\to M_{1,\infty}$ (see also \cite[Theorem II.4.4]{KPS}).

Further, we need to recall the important notion of Hardy-Littlewood majorization. Let $A,B\in\mathcal{N}.$ $B$ is said to be majorized by $A$ and written $B\prec\prec A$ if and only if
\begin{equation}\label{major def}
\int_0^t\mu(s,B)ds\leq\int_0^t\mu(s,A)ds,\quad t\geq0.
\end{equation}
We have (see \cite{FK})
\begin{equation}\label{majorization inequalities}
A+B\prec\prec\mu(A)+\mu(B)\prec\prec2\sigma_{1/2}\mu(A+B).
\end{equation}

One of the most widely used ideals in von Neumann algebras is
$$\mathcal{L}_p:= \mathcal{L}_p(\mathcal{N},\tau)=\{A\in \mathcal{N}:\ \|A\|_p:=\tau(|A|^p)^{1/p}<\infty\}, \ p\ge 1,$$
usually called the Schatten-von Neumann ideal of $p$-summable operators. Using Hardy-Littlewood majorization, it is very easy to see (e.g. \cite[Lemma 2.1]{CPS}) that $\mathcal{M}_{1,\infty}\subset \mathcal{L}_p$ for all $p>1.$

A linear functional $\varphi:\mathcal{M}_{1,\infty}\to\mathbb{C}$ is said to be {\it symmetric} if $\varphi(B)=\varphi(A)$ for every
positive $A,B\in\mathcal{M}_{1,\infty}$ such that $\mu(B)=\mu(A).$ A linear functional $\varphi:\mathcal{M}_{1,\infty}\to\mathbb{C}$ is said to be {\it fully symmetric} if $\varphi(B)\leq\varphi(A)$ for all $A,B\in\mathcal{M}^+_{1,\infty}$ such that $B\prec\prec A$ \cite{DPSSS1,DPSSS2,DPSS}. Every fully symmetric functional is symmetric and bounded. The converse fails \cite{KS}.

A positive normalised linear functional $\gamma:L_{\infty}(0,\infty)\to\mathbb{R}$ is called a {\it generalised limit} if $\gamma(z)=0$ for every $z\in L_{\infty}(0,\infty)$ such that $\lim_{t\to\infty}z(t)=0.$ A linear functional $\gamma:L_{\infty}(0,\infty)\to\mathbb{R}$ is called {\it dilation invariant} if $\gamma(\sigma_sz)=\gamma(z)$ for every $z\in L_{\infty}(0,\infty)$ and every $s>0.$

Let $S\subseteq B(H).$ We denote by $S^+$ the set of all positive operators from $S.$

Let $\omega:L_{\infty}(0,\infty)\to\mathbb{R}$ be a dilation invariant generalised limit. Define a functional $\tau_{\omega}$ on $\mathcal{M}^+_{1,\infty}$ by the formula
$$\tau_{\omega}(A)=\omega(\frac1{\log(1+t)}\int_0^t\mu(s,A)ds).$$
The functional $\tau_{\omega}$ is additive and unitarily invariant on $\mathcal{M}^+_{1,\infty}.$ Thus, $\tau_{\omega}$ extends to a
fully symmetric functional on $\mathcal{M}_{1,\infty}.$ One usually refers to it as to a {\it Dixmier trace}. We refer the reader to \cite{Dixmier, Connes, CPS, CaS, CRSS, KSS} for details.

Further, we use the following properties of Dixmier traces. Let $A\in\mathcal{M}_{1,\infty}$ and let $B\in\mathcal{N}.$ We have
(see \cite{Connes, CPS})
\begin{equation}\label{dixmier trace pr1}
\tau_{\omega}(AB)=\tau_{\omega}(BA).
\end{equation}
Suppose that $B>0.$ It follows from \eqref{dixmier trace pr1} that
\begin{equation}\label{dixmier trace pr2}
\tau_{\omega}(AB)=\tau_{\omega}(B^{1/2}AB^{1/2}).
\end{equation}

Suppose that the trace $\tau$ on the von Neumann algebra $\mathcal{N}$ is infinite and the algebra $\mathcal{N}$ is either diffuse (that is with no minimal projections) or else is $B(H).$ Given any finite sequence $\{A_n\}$ of operators, we can construct a sequence of operators $\{B_n\}$ such that $\mu(A_n)=\mu(B_n)$ for all $n$'s and $B_nB_m=0$ for all $n\neq m.$ Further, we refer to any such sequence $\{B_n\}$ as a ''sequence of disjoint copies of $\{A_n\}$''.

Cesaro operator $M$ is defined on $L_{\infty}(0,\infty)$ by the formula
$$(Mx)(t)=\frac1{\log(t)}\int_1^tx(s)\frac{ds}{s},\quad t\in(0,\infty).$$

\section{Preliminary important results}\label{prel}

In this section, for the reader's convenience, we collect a number of key known results, which will be used throughout this paper.

The following important theorem is proved in \cite[Theorem 11]{KSS} for general Marcinkiewicz spaces.

\begin{thm}\label{fs is dixmier} Every fully symmetric functional on $\mathcal{M}_{1,\infty}$ is a Dixmier trace.
\end{thm}

The following theorem is an analog of Lidskii formula (see \cite{Simon}) for Dixmier traces. It is proved in \cite[Theorem 33]{SSZ} for a large subclass of Marcinkiewicz spaces which contains $\mathcal{M}_{1,\infty}.$

\begin{thm}\label{lidskii formula} Let $A\in\mathcal{M}_{1,\infty}$ and let $\tau_{\omega}$ be an arbitrary Dixmier trace on $\mathcal{M}_{1,\infty}.$ We have
$$\tau_{\omega}(A)=\omega\left(\frac1{\log(t)}\sum_{|\lambda|>\log(t)/t,\lambda\in\sigma(A)}\lambda\right).$$
\end{thm}

The following $\omega$-variant of the classical Karamata theorem is established in \cite{CPS}.

\begin{thm}\label{karamata} Let $\beta$ be a continuous increasing function. Set
$$h(t)=\int_0^{\infty}e^{-(u/t)^q}d\beta(u).$$
We have
$$\omega(\frac{h(t)}{t})=\Gamma(1+\frac1q)\omega(\frac{\beta(t)}{t})$$
for any dilation invariant generalised limit $\omega.$
\end{thm}

Consider the ideal $\mathcal{K}_{\mathcal {N}}$ of $\tau$-compact operators in $\mathcal {N}$ (that is the norm closed ideal generated by the projections $E\in \mathcal {N}$ with $\tau(E) < \infty$). The following result is not new (see \cite[Chapter II, Lemma 3.4]{GK}). We present a short proof for convenience of the reader.

\begin{thm}\label{t equals majorization} Let $A,B\in \mathcal {N}$ be positive $\tau-$compact operators. We have $B\prec\prec A$ if and only if
\begin{equation}\label{other majorization}
\tau((B-t)e^B(t,\infty))\leq\tau((A-t)e^A(t,\infty)),\quad \forall t>0.
\end{equation}
\end{thm}
\begin{proof} Fix $t>0.$ It follows from the definition of generalised singular value function that $\mu(Ae^A(t,\infty))=\mu(A)\chi_{[0,d_A(t)]}.$ Applying \cite[Proposition 2.7]{FK} to the operator $Ae^A(t,\infty),$ we have
$$\tau(Ae^A(t,\infty))=\int_0^{d_A(t)}\mu(s,A)ds,$$
and hence
\begin{equation}\label{t repres}
\tau((A-t)e^A(t,\infty))=\int_0^{d_A(t)}(\mu(s,A)-t)ds.
\end{equation}
The function
$$u\to\int_0^u(\mu(s,A)-t)ds$$
attains its maximum at $u=d_A(t).$

If $B\prec\prec A,$ then
$$\int_0^{d_B(t)}(\mu(s,B)-t)ds\leq\int_0^{d_B(t)}(\mu(s,A)-t)ds\leq\int_0^{d_A(t)}(\mu(s,A)-t)ds.$$
Inequality \eqref{other majorization} follows now from \eqref{t repres}.

Suppose now that \eqref{other majorization} holds. Fix $u>0$ and set $t=\mu(u,A).$ It follows that
$$\int_0^u(\mu(s,B)-t)ds\leq\int_0^{d_B(t)}(\mu(s,B)-t)ds=\tau((B-t)e^B(t,\infty))\leq$$
$$\leq\tau((A-t)e^A(t,\infty))=\int_0^u(\mu(s,A)-t)ds.$$
Hence,
$$\int_0^u\mu(s,B)ds\leq\int_0^u\mu(s,A)ds.$$
Since $u$ is arbitrary, we have $B\prec\prec A.$
\end{proof}

\section{$\zeta-$function formulae}\label{zeta}

We begin by showing that the functionals given in \eqref{gamma-zeta-def} are well defined on $\mathcal{M}^+_{1,\infty}.$

\begin{lem}\label{zeta correctness} If $\gamma:L_{\infty}(0,\infty)\to\mathbb{R}$ is a generalised limit, then $\zeta_{\gamma}(A)<\infty$ and $\zeta_{\gamma,B}(A)<\infty$ for any $A\in\mathcal{M}^+_{1,\infty}.$
\end{lem}
\begin{proof} It is clear that $\mu(s,A)\prec\prec(1+s)^{-1}\|A\|_{1,\infty}.$ Therefore,
$$\tau(A^{1+1/t})\leq\|A\|_{1,\infty}^{1+1/t}\int_0^{\infty}\frac{dt}{(1+s)^{1+1/t}}=t\|A\|_{1,\infty}^{1+1/t}.$$
Hence, $\zeta_{\gamma}(A)\leq\|A\|_{1,\infty}.$ It follows from
$$\tau(A^{1+1/t}B)\leq\|B\|\tau(A^{1+1/t})$$
that $\zeta_{\gamma,B}(A)\leq \|B\|\zeta_{\gamma}(A).$
\end{proof}

\begin{rem}\label{chei remark} Let $x,y\in L_{\infty}(0,\infty).$ For any generalised limit $\gamma$ such that $\gamma(|x-1|)=0,$ we have $\gamma(xy)=\gamma(y).$ Indeed, $|\gamma(xy-y)|\leq\gamma(|x-1|)\|y\|=0.$
\end{rem}

\begin{lem}\label{elementary inequality} For any $A,C\in \mathcal{M}^+_{1,\infty}$ we have
$$\tau(A^{1+s}+C^{1+s})\leq\tau((A+C)^{1+s})\leq2^s\tau(A^{1+s}+C^{1+s}),\quad s>0.$$
\end{lem}
\begin{proof} In the special case when $\mathcal{N}=B(H),$ the first inequality can be found in \cite[(2.9)]{Ko}. In the general case, it follows directly from Proposition 4.6(ii) of \cite{FK} when $f(u)=u^{1+s},$ $u>0.$ The second inequality follows from the same proposition by setting there $a=a^*=b=b^*=2^{-1/2}.$
\end{proof}

Let $A\in \mathcal{M}_{1,\infty}.$ For a functional $\zeta_{\gamma}$ defined on $\mathcal{M}^+_{1,\infty}$ by \eqref{gamma-zeta-def} (see Lemma \ref{zeta correctness}), we set
\begin{equation}\label{zeta-full}
\zeta_{\gamma}(A):=(\zeta_{\gamma}(\Re(A)_+)-\zeta_{\gamma}(\Re(A)_-))+i(\zeta_{\gamma}(\Im(A)_+)-\zeta_{\gamma}(\Im(A)_-)).
\end{equation}

The following theorem shows that functionals $\zeta_{\gamma}$ defined by \eqref{zeta-full} are fully symmetric on $\mathcal{M}_{1,\infty}.$

\begin{thm}\label{zeta linearity} If $\gamma:L_{\infty}(0,\infty)\to\mathbb{R}$ is a generalised limit, then  $\zeta_{\gamma}$ is a fully symmetric linear functional on $\mathcal{M}_{1,\infty}.$
\end{thm}
\begin{proof} To verify that   $\zeta_{\gamma}$ is linear, it is sufficient to check that
$\zeta_{\gamma}(A+C)=\zeta_{\gamma}(A)+\zeta_{\gamma}(C)$ for any $A,C\in \mathcal{M}^+_{1,\infty}.$  It follows from the left hand side inequality of Lemma \ref{elementary inequality} that
$$\zeta_{\gamma}(A+C)\geq\zeta_{\gamma}(A)+\zeta_{\gamma}(C).$$
Noting that $\gamma(|2^{1/t}-1|)=0,$ it follows from the right hand side inequality of Lemma \ref{elementary inequality} and Remark \ref{chei remark} that
$$\zeta_{\gamma}(A+C)\leq\zeta_{\gamma}(A)+\zeta_{\gamma}(C).$$

Therefore, we have
$$\zeta_{\gamma}(A+C)=\zeta_{\gamma}(A)+\zeta_{\gamma}(C).$$
The homogeneity of $\zeta_{\gamma}$ follows from Remark \ref{chei remark}.
Finally, if $0\leq C\prec\prec A\in\mathcal{M}^+_{1,\infty},$ then $C,A\in\mathcal{L}_{1+s}$ and $\tau(C^{1+s})\leq\tau(A^{1+s}).$ Hence,
$\frac1t\tau(C^{1+1/t})\leq\frac1t\tau(A^{1+1/t})$ and so $\zeta_{\gamma}(C)\leq \zeta_{\gamma}(A).$
\end{proof}

Let $B\in \mathcal{N}.$ We extend the functional $\zeta_{\gamma,B}$ on $\mathcal{M}_{1,\infty},$ similarly to \eqref{zeta-full}. Observe that
$$\zeta_{\gamma,B_1+B_2}(A)=\zeta_{\gamma,B_1}(A)+\zeta_{\gamma,B_2}(A),\ B_1, B_2\in \mathcal{N},\ A\in \mathcal{M}_{1,\infty}.$$

\begin{lem}\label{continuity lemma} If $A\in\mathcal{M}_{1,\infty}$ and $B_n\to B$ in $\mathcal{N},$ then
$$\zeta_{\gamma,B_n}(A)\to\zeta_{\gamma,B}(A).$$
\end{lem}
\begin{proof} It is sufficient to prove the assertion for $A\in\mathcal{M}^+_{1,\infty}.$ Since
$$|\tau(A^{1+s}B)-\tau(A^{1+s}B_n)|\leq\tau(A^{1+s})\|B-B_n\|,$$
we obtain
$$|\zeta_{\gamma,B}(A)-\zeta_{\gamma,B_n}(A)|\leq\zeta_{\gamma}(A)\|B-B_n\|.$$
\end{proof}

The following lemma follows immediately from \cite[Lemma 3.3]{CPS}.

\begin{lem}\label{CPS lemma} Let $A,B\in B^+(H)$ and let $s>0.$ We have
\begin{enumerate}[i)]
\item $(B^{1/2}AB^{1/2})^{1+s}\leq B^{1/2}A^{1+s}B^{1/2}$ if $0\leq B\leq 1.$
\item $(B^{1/2}AB^{1/2})^{1+s}\geq B^{1/2}A^{1+s}B^{1/2}$ if $B\geq 1.$
\end{enumerate}
\end{lem}

The result below significantly strengthens \cite[Proposition 3.6]{CPS} by removing all extra assumptions on the generalised limit $\gamma.$

\begin{prop}\label{main reduction lemma} If  $\gamma:L_{\infty}(0,\infty)\to\mathbb{R}$ is a generalised limit, then
$$\zeta_{\gamma,B}(A)=\zeta_{\gamma}(B^{1/2}AB^{1/2}),\ \forall A\in\mathcal{M}_{1,\infty},\ B\in \mathcal{N}^+.$$
\end{prop}
\begin{proof} It is sufficient to prove the assertion for $A\in\mathcal{M}^+_{1,\infty}.$
Suppose first that there are constants $0<m\leq M<\infty$ such that $m\leq B\leq M.$ Applying Lemma \ref{CPS lemma} to the operators $A$ and $M^{-1}B$ (respectively, $m^{-1}B$), we have
$$m^sB^{1/2}A^{1+s}B^{1/2}\leq (B^{1/2}AB^{1/2})^{1+s}\leq M^sB^{1/2}A^{1+s}B^{1/2}.$$
Therefore,
$$\frac1tm^{1/t}\tau(A^{1+1/t}B)\leq\frac1t\tau((B^{1/2}AB^{1/2})^{1+1/t})\leq\frac1t M^{1/t}\tau(A^{1+1/t}B).$$
Since $\gamma(|m^{1/t}-1|)=0$ and $\gamma(|M^{1/t}-1|)=0,$  it follows from Remark \ref{chei remark} that $\zeta_{\gamma,B}(A)=\zeta_{\gamma}(B^{1/2}AB^{1/2}).$

For an arbitrary $B\in \mathcal{N}^+,$ we set $B_n:=Be^B(1/n,\infty)+1/ne^B[0,1/n],$ $n\ge 1.$ From the first part of the proof, we have
$$\zeta_{\gamma,B_n}(A)=\zeta_{\gamma}(B_n^{1/2}AB_n^{1/2}).$$
Since $B_n^{1/2}AB_n^{1/2}\to B^{1/2}AB^{1/2}$ in $\mathcal{M}_{1,\infty},$ we have by Theorem \ref{zeta linearity}
$$\zeta_{\gamma}(B_n^{1/2}AB_n^{1/2})\to\zeta_{\gamma}(B^{1/2}AB^{1/2}).$$
On the other hand, by Lemma \ref{continuity lemma} we have $\zeta_{\gamma,B_n}(A)\to\zeta_{\gamma,B}(A).$
\end{proof}

The following is our main result on the $\zeta-$function.

\begin{thm}\label{zeta reduction} If  $\gamma:L_{\infty}(0,\infty)\to\mathbb{R}$ is a generalised limit, then
$$\zeta_{\gamma,B}(A)=\zeta_{\gamma}(AB),\ \forall A\in\mathcal{M}_{1,\infty},\ B\in \mathcal{N}.$$
\end{thm}
\begin{proof} It is sufficient to prove the assertion for $B\in \mathcal{N}^+.$ By Theorems \ref{zeta linearity}  and \ref{fs is dixmier}, we know that $\zeta_{\gamma}$ is a Dixmier trace on $\mathcal{M}_{1,\infty}.$ Hence, by \eqref{dixmier trace pr2}, we have $\zeta_{\gamma}(B^{1/2}AB^{1/2})=\zeta_{\gamma}(AB).$ The assertion follows now from Proposition \ref{main reduction lemma}.
\end{proof}

Our remaining objective in this section is to provide
strengthening of several formulae linking Dixmier traces and
$\zeta$-functions from \cite{CPS, CRSS}.

\begin{lem}\label{convexity lemma} Let $A\in\mathcal{M}_{1,\infty}^+.$ The mapping $s\to s^{-1}\zeta_{\gamma\circ\sigma_s}(A)$ is convex and, therefore, continuous.
\end{lem}
\begin{proof} For all $t,s>0,$ we have
$$s^{-1}\sigma_s(\frac1t\tau(A^{1+1/t}))=\frac1t\tau(A^{1+s/t}).$$
Therefore, for every $s>0$
$$s^{-1}\zeta_{\gamma\circ\sigma_s}=\gamma(\frac1t\tau(A^{1+s/t})).$$
Let $\lambda_i>0$ and let $\lambda_1+\lambda_2=1.$ Since the mapping $t\to a^{1+t}$ is convex for every $a>0,$ it follows from the spectral theorem that the map $s\to A^s$ is also convex. Therefore, for all positive real numbers $s_1,s_2 $ and $t,$ we have
$$A^{1+(\lambda_1s_1+\lambda_2s_2)/t}\leq\lambda_1A^{1+s_1/t}+\lambda_2A^{1+s_2/t}.$$
The assertion follows immediately.
\end{proof}

Let $\gamma$ be a generalised limit on $L_{\infty}(0,\infty).$ Below, we will formally apply the notation $\zeta_{\gamma,B}(A)$ introduced in \eqref{gamma-zeta-def} to some unbounded positive operators $B$ on $H.$

\begin{lem}\label{simple technical estimate} Let $A\in\mathcal{N}$ be a positive $\tau-$compact operator and let $B\geq1$ be an unbounded operator commuting with $A.$ If (the closure of) the product $AB\in\mathcal{M}_{1,\infty}$ and $AB^n\in\mathcal{N}$ for every $n\in\mathbb{N},$ then $\zeta_{\gamma}(AB)=\zeta_{\gamma,B}(A).$
\end{lem}
\begin{proof} It follows from $AB=BA$ and $B\geq 1$ that $A^{1+s}B\leq(AB)^{1+s}.$ The inequality $\zeta_{\gamma,B}(A)\leq\zeta_{\gamma}(AB)$ follows immediately.

Set $c_n:=\|AB^{2n}\|,$ $n\ge 1$ and observe that $BA^{1/2n}\leq c_n^{1/2n}.$ Setting $B_n=Be^A[0,c_n^{-1}],$ we obtain
\begin{equation}\label{b leq a1n}
B_nA^{1/n}=BA^{1/2n}\cdot A^{1/2n}e^A[0,c_n^{-1}]\leq(c_nA)^{1/2n}e^A[0,c_n^{-1}]\leq 1.
\end{equation}

It follows from \eqref{b leq a1n} that $A^{1+1/t}B_n\geq (AB_n)^{1+n/t(n-1)}.$ Thus,
$$\gamma(\frac1t\tau(A^{1+1/t}B_n))\geq\gamma(\frac1t\tau((AB_n)^{1+n/t(n-1)}))=\frac{n-1}{n}\zeta_{\gamma\circ\sigma_{n/(n-1)}}(AB_n).$$

Since $A$ is $\tau-$compact, then $B-B_n$ is bounded operator with finite support. Due to the linearity with respect to $B,$ we have
$$\zeta_{\gamma,B}(A)=\zeta_{\gamma,B_n}(A)\geq\frac{n-1}{n}\zeta_{\gamma\circ\sigma_{n/(n-1)}}(AB_n)=\frac{n-1}{n}\zeta_{\gamma\circ\sigma_{n/(n-1)}}(AB).$$
The assertion follows now from Lemma \ref{convexity lemma}.
\end{proof}

The following result is mainly known (see \cite{CPS, CRSS}). Our proof is however much simpler than the arguments used there.

\begin{thm}\label{zeta dix}  If $\omega$ is a dilation invariant generalised limit such that the generalised limit $\omega\circ\log$ is still dilation invariant, then $\tau_{\omega}=\zeta_{\omega\circ\log}.$
\end{thm}
\begin{proof} It is sufficient to verify the equality $\tau_{\omega}=\zeta_{\omega\circ\log}$ on positive operators $A\in \mathcal{M}_{1,\infty}^+$ such that $A\leq e^{-1}.$ Define a continuously increasing function $\beta:(0,\infty)\to(0,\infty)$ by
$$\beta(u):=-\int_{ue^{-u}}^{\infty}\lambda dd_A(\lambda).$$
Let $h$ be as in Theorem \ref{karamata} as applied to the above
$\beta.$ Define an operator $B\geq 1$ by the formula $A=Be^{-B}$
and set $C=e^{-B}.$ We have
\begin{equation}\label{h(t)/t}
h(t)=\int_0^{\infty}e^{-u/t}d\beta(u)=-\int_0^{\infty}e^{-u(1+1/t)}udd_A(ue^{-u})\stackrel{\eqref{f formula}}{=}\tau(C^{1+1/t}B).
\end{equation}

The conditions of Lemma \ref{simple technical estimate} are valid for $B$ and $C.$ Indeed, $B$ commutes with $C,$ $BC=A\in\mathcal{M}_{1,\infty}$ and $B^ne^{-B}\in\mathcal{N}$ for every $n\in\mathbb{N}.$ By Lemma \ref{simple technical estimate}, we have
$$\zeta_{\omega\circ\log}(A)=\zeta_{\omega\circ\log,B}(C)=(\omega\circ\log)(\frac{h(t)}{t}).$$
By Theorem \ref{lidskii formula}, we have
\begin{equation}\label{beta(t)/t} \tau_{\omega}(A)=\omega(\frac{-1}{\log(t)}\int_{\log(t)/t}^{\infty}\lambda dd_A(\lambda))=(\omega\circ\log)(\frac{\beta(t)}{t}).
\end{equation}
We can now conclude
$$\zeta_{\omega\circ\log}(A)\stackrel{\mathrm{ \eqref{h(t)/t}}}{=}(\omega\circ\log)(\frac{h(t)}{t})\stackrel{\mathrm{(Thm \, \ref{karamata})}}{=}(\omega\circ\log)(\frac{\beta(t)}{t})\stackrel{\mathrm{( \ref{beta(t)/t})}}{=}\tau_{\omega}(A).$$
\end{proof}

The following corollary strengthens and extends the results of
\cite[Theorem 4.11]{CRSS} and \cite[Theorem 3.8]{CPS}.
It follows immediately from Theorems \ref{zeta dix} and \ref{zeta reduction}.

\begin{cor} If $\omega$ is a dilation invariant generalised limit such that the generalised limit
$\omega\circ\log$ is still dilation invariant, then
$$\tau_{\omega}(AB)=(\omega\circ\log)(\frac1t\tau(A^{1+1/t}B)),\ \forall A\in\mathcal{M}^+_{1,\infty},\
B\in \mathcal{N}.$$
\end{cor}

\section{The linearity criterion for functionals $\xi_\gamma$}\label{heatk lin}

In this section we focus on functionals $\xi_{\gamma}(\cdot)$ defined in \eqref{modified-gamma-hk-q-def}. It was implicitly proved in \cite[Theorem 5.2]{CRSS} that
$$M\left(t\to \frac1t\tau(\exp(-(tA)^{-q}))\right)\in L_{\infty}(0,\infty),\ \forall A\in\mathcal{M}^+_{1,\infty}$$
and therefore,
\begin{equation}\label{xi-repeat}
\xi_{\gamma}(A):=(\gamma\circ M)\left(t\to \frac1t\tau(\exp(-(tA)^{-q}))\right)
\end{equation}
is finite for every $A\in \mathcal{M}^{+}_{1,\infty}$ and every generalised limit $\gamma$ on $L_{\infty}(0,\infty).$ We note, in passing that a stronger result than \cite[Theorem 5.2]{CRSS} is established in Theorem \ref{mhk bounded} below. Let $A\in \mathcal{M}_{1,\infty}.$ For a functional $\xi_{\gamma},$ we set
\begin{equation}\label{xi-full}
\xi_{\gamma}(A):=(\xi_{\gamma}(\Re(A)_+)-\xi_{\gamma}(\Re(A)_-))+i(\xi_{\gamma}(\Im(A)_+)-\xi_{\gamma}(\Im(A)_-)).
\end{equation}
It is probably a difficult task to describe the set of all generalised limits $\gamma$ for which \eqref{xi-full} yields a linear functional $\xi_{\gamma}.$ However, the class of linear functionals $\xi_{\gamma}$ is an easier object. Below in Proposition
\ref{coincidence}, we show that the sets of linear functionals $\{\xi_\gamma:\ \gamma\ \text{is\ a generalised\ limit}\}$ and linear functionals $\{\xi_\omega:\ \omega\ \text{is\ a\ dilation\ invariant\ generalised\ limit}\}$ coincide.

\begin{lem}\label{stupid m remark} For every locally integrable $z$ with $Mz\in L_{\infty}(0,\infty),$ we have
$$(M\circ\sigma_{s^{-1}}-\sigma_{s^{-1}}\circ M)(z)\in C_0^b(0,\infty),\ \forall s>0.$$
Here, $C_0^b(0,\infty)$ is the space of all bounded continuous functions tending to $0$ at $\infty.$
\end{lem}
\begin{proof} Fix $s>0.$ The assertion follows by writing
$$(M\circ\sigma_{s^{-1}}-\sigma_{s^{-1}}\circ M)(z)
=\frac1{\log(t)}\int_{s}^{st}z(u)\frac{du}{u}-\frac1{\log(st)}\int_1^{st}z(u)\frac{du}{u}$$
and noting that the assumption $Mz\in L_{\infty}(0,\infty)$ easily implies that
$$\frac1{\log(st)}\int_1^{st}z(u)\frac{du}{u}-\frac1{\log(t)}\int_1^{st}z(u)\frac{du}{u}\in C_b^0(0,\infty).$$
\end{proof}

\begin{prop}\label{coincidence} Suppose that a generalised limit $\gamma$ on $L_{\infty}(0,\infty)$ is such that $\xi_{\gamma}$ is a linear functional on $\mathcal{M}_{1,\infty}.$ Then, there exists a dilation invariant generalised limit $\omega$ on $L_{\infty}(0,\infty)$ such that $\xi_{\gamma}=\xi_{\omega}.$
\end{prop}
\begin{proof} Fix $s>0$ and observe that
\begin{equation}\label{bl1}
\left(t\to
\frac1t\tau(\exp(-(tsA)^{-q}))\right)=s\sigma_{s^{-1}}\left(t\to
\frac1t\tau(\exp(-(tA))^{-q})\right).
\end{equation}
Therefore,
$$\xi_{\gamma}(sA)=s(\gamma\circ M\circ\sigma_{s^{-1}})(\frac1t\tau(\exp(-(tA)^{-q}))).$$
By the assumption, we have $\xi_{\gamma}(sA)=s\xi_{\gamma}(A)$ and
appealing to Lemma \ref{stupid m remark}, we obtain
\begin{equation}\label{xi as dilation}
\xi_{\gamma}(A)=(\gamma\circ\sigma_{s^{-1}}\circ M)(\frac1t\tau(\exp(-(tA)^{-q}))),\ \forall s>0.
\end{equation}
Let $E$ be the linear span of the functions
$$t\to M(\frac1t\tau(\exp(-(tA)^{-q}))),\ A\in \mathcal{M}^+_{1,\infty}$$
and let $F:=E+C_0^b(0,\infty).$ We claim that the space $F$ is dilation invariant. Indeed, it follows from Lemma \ref{stupid m remark} and \eqref{bl1} that every function
$$\sigma_{s^{-1}}\left(t\to M(\frac1t\tau(\exp(-(tA))^{-q}))\right)$$
belongs to the set
$$s^{-1}\left(t\to M(\frac1t\tau(\exp(-(tsA)^{-q})))\right)+C_0^b(0,\infty).$$
It follows from \eqref{xi as dilation} that $\gamma\circ\sigma_{s^{-1}}=\gamma$ on $F.$ By the invariant form of the Hahn-Banach theorem (see \cite[p. 157]{Edwards}) applied to the group of dilations $\{\sigma_s\}_{s>0},$ we see that $\gamma|_F$ can be extended to a dilation invariant generalised limit $\omega$ on $L_\infty(0,\infty).$
\end{proof}

The following lemma can be found in \cite{SSZ}. We present a shorter proof for convenience of the reader.

\begin{lem}\label{xiomega repr0} If $\omega$ is a dilation invariant generalised limit on $L_{\infty}(0,\infty),$ then
\begin{equation}\label{hk karamata equality}
\xi_{\omega}(A)=\Gamma(1+\frac1q)(\omega\circ
M)(\frac1td_A(\frac1t)),\ \forall A\in\mathcal{M}^+_{1,\infty}.
\end{equation}
\end{lem}
\begin{proof} It follows from \eqref{f formula} that
\begin{equation}\label{add}
\tau(\exp(-(tA)^{-q}))=\int_0^{\infty}e^{-(u/t)^q}dd_A(\frac1u).
\end{equation}
Setting $\beta(u)=d_A(1/u),$ multiplying both sides of \eqref{add} by $1/t$ and applying Theorem \ref{karamata} to $\omega\circ M$ (which is dilation invariant, see \cite{Connes}), we obtain \eqref{hk karamata equality}.
\end{proof}

\begin{lem}\label{xiomega repr} Let $A\in\mathcal{M}^+_{1,\infty}$ and let $\omega$ be a dilation invariant generalised limit on $L_{\infty}(0,\infty).$ We have
\begin{equation}\label{hk karamata equality2}
\xi_{\omega}(A)=\Gamma(1+\frac1q)\omega(\frac1{\log(1+t)}\tau((A-\frac1t)e^A(\frac1t,\infty))).
\end{equation}
\end{lem}
\begin{proof} In view of Lemma \ref{xiomega repr0}, it is sufficient to show that right hand sides of \eqref{hk karamata equality} and \eqref{hk karamata equality2} coincide. This easily follows from the following computation, where we use integration by parts
$$M(\frac1td_A(\frac1t))=\frac1{\log(t)}\int_1^td_A(\frac1s)\frac{ds}{s^2}=\frac1{\log(t)}\int_{1/t}^1d_A(u)du=$$
$$=\frac1{\log(t)}ud_A(u)|_{1/t}^1-\frac1{\log(t)}\int_{1/t}^1udd_A(u)=\frac1{\log(t)}\tau((A-\frac1t)e^A(\frac1t,\infty))+o(1).$$
\end{proof}

\begin{lem}\label{full symmetricity for hk} Let $\omega$ be a dilation invariant generalised limit on $L_{\infty}(0,\infty)$ and let $A,B\in\mathcal{M}^+_{1,\infty}$ be such that $B\prec\prec A.$ We have $\xi_{\omega}(B)\leq\xi_{\omega}(A).$
\end{lem}
\begin{proof} The assertion follows from Lemma \ref{xiomega repr} and Theorem \ref{t equals majorization}.
\end{proof}

The following is the main result of this section.

\begin{thm}\label{xiomega linearity} For any dilation invariant generalised limit $\omega$ on $L_{\infty}(0,\infty),$ the functional $\xi_{\omega}$ given by \eqref{xi-full} is linear and fully symmetric on $\mathcal{M}_{1,\infty}.$
\end{thm}
\begin{proof} The assertion follows from Lemma \ref{full symmetricity for hk} provided we have shown that
\begin{equation}\label{once again}
\xi_{\omega}(A+B)=\xi_{\omega}(A)+\xi_{\omega}(B),\ \forall A,B\in \mathcal{M}^+_{1,\infty}.
\end{equation}
To this end, we observe first that since $\omega$ and $\omega\circ M$ are dilation invariant, it follows from Lemma \ref{full symmetricity for hk} and \eqref{majorization inequalities} that
$$
\xi_{\omega}(A+B)=\xi_{\omega}(\mu(A)+\mu(B)),\ \forall A,B\in \mathcal{M}^+_{1,\infty}.
$$
Now, let $C$ and $D$ be disjoint copies of $A$ and $B$ (see Section \ref{definot}). Thus, we have
$$\xi_{\omega}(C+D)=\xi_{\omega}(\mu(C)+\mu(D))=\xi_{\omega}(\mu(A)+\mu(B))=\xi_{\omega}(A+B).$$
However, the equality
$$\xi_{\omega}(C+D)=\xi_{\omega}(C)+\xi_{\omega}(D)$$
for positive operators $C$ and $D$ such that $CD=0$ follows immediately from the definition \eqref{xi-repeat}. Since the equalities $\xi_{\omega}(A)=\xi_{\omega}(C),$  $\xi_{\omega}(B)=\xi_{\omega}(D)$ are obvious, we arrive at \eqref{once again}.
\end{proof}

\section{Every fully symmetric functional has form $\xi_\omega$}\label{heatk all}

It follows from Theorem \ref{xiomega linearity} and Theorem \ref{fs is dixmier}, that the functional $\xi_{\omega}$ is a fully symmetric functional on $\mathcal{M}_{1,\infty}$ whenever $\omega$ is a dilation invariant generalised limit $\omega$ on $L_{\infty}(0,\infty).$ In this section, we show the converse.

Define a (non-linear) operator $T:\mathcal{M}^+_{1,\infty}\to L_{\infty}(0,\infty)$ by the formula
\begin{equation}\label{TA}(TA)(t)=\frac1{\log(1+t)}\tau((A-\frac1t)e^A(\frac1t,\infty)),\ t>0.
\end{equation}

We need some properties of the operator $T.$ Firstly, we show that it is additive on certain pairs of $A,B\in \mathcal{M}^+_{1,\infty}.$

\begin{lem}\label{disjoint linearity} Let $A,B\in \mathcal{M}^+_{1,\infty}$ be such that $AB=BA=0.$ It follows that $T(A+B)=TA+TB.$
\end{lem}
\begin{proof} It follows immediately from the assumption that
$$(A+B-\frac1t)e^{A+B}(\frac1t,\infty)=(A-\frac1t)e^A(\frac1t,\infty)+(B-\frac1t)e^B(\frac1t,\infty).$$
\end{proof}

Next, we explain the connection of the operator $T$ with fully symmetric functionals on $\mathcal{M}_{1,\infty}.$

\begin{lem}\label{fs ab estimate} Let the operators $A,B\in\mathcal{M}^+_{1,\infty}$ be such that $TB\leq TA.$ For every fully symmetric functional $\varphi$  on $\mathcal{M}_{1,\infty},$ we have $\varphi(B)\leq\varphi(A).$
\end{lem}
\begin{proof} It follows immediately from the definition \eqref{TA} that
$$\tau((B-\frac1t)e^B(\frac1t,\infty))\leq\tau((A-\frac1t)e^A(\frac1t,\infty)),\quad \forall t>0.$$
Applying Theorem \ref{t equals majorization} we obtain $B\prec\prec A$ and so $\varphi(B)\leq\varphi(A).$
\end{proof}

\begin{lem}\label{razn1} Let $A,B\in\mathcal{M}^+_{1,\infty}.$ For every fully symmetric functional $\varphi$  on $\mathcal{M}_{1,\infty},$ we have
$$\varphi(B)-\varphi(A)\leq\|\varphi\|_{\mathcal{M}_{1,\infty}^*}\limsup_{t\to\infty}(TB-TA)(t).$$
\end{lem}
\begin{proof} Without loss of generality, $\|\varphi\|_{\mathcal{M}_{1,\infty}^*}=1.$ Denote the right hand side by $c$ and suppose that $c\geq0$ (the case when $c<0$ is treated similarly). Fix $\varepsilon>0.$ We have $(TB-TA)(t)\leq c+\varepsilon$ for all sufficiently large $t.$ Let $C$ be an operator with $\mu(t,C)=(c+2\varepsilon)/(1+t).$ We have $TB\leq TA+TC$ for all sufficiently large $t.$ Let $A_1$ and $C_1$ be disjoint copies of $A$ and $C,$ respectively. It follows from Lemma \ref{disjoint linearity} that $TB(t)\leq T(A_1+C_1)(t)$ for all sufficiently large $t.$ Choose $0<\delta$ small enough to guarantee $TB_1(t)\leq T(A_1+C_1)(t)$ for all $t>0,$ where $B_1:=\min\{B,\delta\}.$ By Corollary \ref{fs ab estimate}, we have $\varphi(B_1)\leq\varphi(A_1)+\varphi(C_1),$ or equivalently  $\varphi(B)\leq\varphi(A)+c+2\varepsilon.$ Since $\varepsilon$ is arbitrarily small, we are done.
\end{proof}

\begin{lem}\label{ro est} Let $A_1,\cdots,A_n\in\mathcal{M}^+_{1,\infty}$ and let $\lambda_1,\cdots,\lambda_n\in\mathbb{R}$ for some $n\ge 1.$ For every fully symmetric functional $\varphi$  on $\mathcal{M}_{1,\infty}$ we have
\begin{equation}\label{ro estimate}
\sum_{k=1}^n\lambda_k\varphi(A_k)\leq\limsup_{t\to\infty}\sum_{k=1}^n\lambda_k(TA_k)(t).
\end{equation}
\end{lem}
\begin{proof} Both sides of the inequality \eqref{ro estimate} depend continuously on the $\lambda_k$'s.  Without loss of generality, we may assume that all $\lambda_k\in\mathbb{Q}.$ Multiplying both sides by the common denominator, we may assume that all $\lambda_k\in\mathbb{Z}.$ Writing
$$\lambda_kA_k=\sum_{k=1}^{|\lambda_k|}{\rm sgn}(\lambda_k)A_k$$
we see that it is sufficient to prove \eqref{ro estimate} only for the case when $\lambda_k=\pm 1$ for every $k.$

Let $\{B_k\}$ be a disjoint copy sequence of $\{A_k\}.$ Both sides of the inequality \eqref{ro estimate} do not change if we replace $A_k$ with $B_k.$ Without loss of generality, the operators $A_kA_j=0,$ $k\neq j.$ By Lemma \ref{razn1} we have
$$\sum_{k=1}^n\lambda_k\varphi(A_k)=\varphi(\sum_{\lambda_k=1}A_k)-\varphi(\sum_{\lambda_k=-1}A_k)\leq$$
$$\leq\limsup_{t\to\infty}(T(\sum_{\lambda_k=1}A_k)-T(\sum_{\lambda_k=-1}A_k))(t).$$
Since $A_kA_j=0$ for all $k\neq j,$ we have by Lemma \ref{disjoint linearity} that
$$T(\sum_{\lambda_k=1}A_k)-T(\sum_{\lambda_k=-1}A_k)=\sum_{k=1}^n\lambda_kTA_k$$
and the assertion follows.
\end{proof}

\begin{lem} Let $E$ be the linear span of $T\mathcal{M}^+_{1,\infty}$ and $C_0^b(0,\infty).$ For every $s>0$ we have $\sigma_sE=E.$
\end{lem}
\begin{proof} It follows from the definition \eqref{TA} that for every $s>0,$ we have
\begin{equation}\label{sigma e}
\sigma_sTA\in sT(s^{-1}A)+C_0^b(0,\infty),\ \forall A\in \mathcal{M}^+_{1,\infty}.
\end{equation}
\end{proof}

Let $\varphi$ be a normalised fully symmetric functional on $\mathcal{M}_{1,\infty}.$ We need the following linear functional  on $E.$

\begin{defi}\label{ro defi} For every $z\in E$ such that
$$z\in\sum_{k=1}^n\lambda_kTA_k+C_0^{\infty}(0,\infty)$$
we set
$$\rho(z)=\sum_{k=1}^n\lambda_k\varphi(A_k).$$
\end{defi}

That $\rho$ is well-defined is proved below.

\begin{lem}\label{rho estimate} The linear functional $\rho:E\to\mathbb{R}$ is well-defined. For every $z\in E,$ we have
$$\rho(z)\leq\limsup_{t\to\infty}z(t).$$
\end{lem}
\begin{proof} Let $z\in E$ be such that
$$z\in\sum_{k=1}^n\lambda_kTA_k+C_0^b(0,\infty),\quad z\in\sum_{k=1}^m\mu_kTB_k+C_0^b(0,\infty).$$
We have
$$\sum_{k=1}^n\lambda_kTA_k-\sum_{k=1}^m\mu_kTB_k\in C_0^b(0,\infty).$$
It follows from Lemma \ref{ro est} that
$$\sum_{k=1}^n\lambda_k\varphi(A_k)=\sum_{k=1}^m\mu_k\varphi(B_k),$$
so that $\rho$ is well-defined.

The second assertion directly follows from Lemma \ref{ro est}.
\end{proof}

\begin{lem}\label{main fs lemma} Let $\varphi$ be a normalised fully symmetric functional on $\mathcal{M}_{1,\infty}.$ There exists a dilation invariant generalised limit $\omega$ on $L_{\infty}(0,\infty)$ such that $\varphi(A)=\omega(TA)$ for every $A\in\mathcal{M}^+_{1,\infty}.$
\end{lem}
\begin{proof} For every $A\in\mathcal{M}^+_{1,\infty},$ we have
$$\rho(\sigma_sTA)\stackrel{\eqref{sigma e}}{=}\rho(sT(s^{-1}A))\stackrel{\text{Def.} \ref{ro defi}}{=}s\varphi(s^{-1}A)=\rho(TA).$$
Therefore, $\rho$ is $\sigma_s-$invariant on $E.$ It follows from Lemma \ref{rho estimate} that
$$\rho(z)\leq\limsup_{t\to\infty}z(t),\quad z\in E.$$
By the invariant form of the Hahn-Banach theorem (see \cite[p. 157]{Edwards}) applied to the group of dilations $\{\sigma_s\}_{s>0},$ we can extend $\rho$ to a dilation invariant generalised limit on $L_{\infty}(0,\infty).$
\end{proof}

The following assertion is the main result of this section. It permits representation of a fully symmetric functional $\varphi$ via heat kernel formulae.

\begin{thm}\label{also main} Let $\varphi$ be a fully symmetric functional on $\mathcal{M}_{1,\infty}.$ There exists dilation invariant
generalised limit $\omega$ on $L_{\infty}(0,\infty)$ such that $\varphi=const\cdot\xi_{\omega}.$
\end{thm}
\begin{proof} It follows from Lemma \ref{main fs lemma} that there exists a dilation invariant generalised limit $\omega$ such that
$$\varphi(A)=\omega(\frac1{\log(1+t)}\tau((A-\frac1t)e^A(\frac1t,\infty))).$$
The assertion follows now from Lemma \ref{xiomega repr}.
\end{proof}

\section{A counterexample}

It is known (see \cite[Theorem 33]{SSZ} and the more general result in Corollary \ref{corollary 51} below) that the equality
$$\xi_{\omega}(A)=\Gamma(1+\frac1q)\tau_{\omega}(A),\quad A\in\mathcal{M}_{1,\infty}^+$$
holds for every $M-$invariant generalised limit $\omega$ on $L_{\infty}(0,\infty)$ (see also earlier results with more restrictive assumptions on $\omega$ in  \cite[Theorem 4.1]{CPS} and \cite[Theorem 5.2]{CRSS}). In view of Theorem \ref{also main} and Theorem \ref{fs is dixmier}, it is quite natural to ask whether the equality above holds for every dilation invariant generalised limit $\omega.$  In this section we prove that this is not the case.

\begin{lem}\label{est1} Let $\omega$ be a dilation invariant generalised limit on $L_{\infty}(0,\infty).$ For every $s>1,$ we have
\begin{equation}\label{eqqq1}
\omega(\sum_k\chi_{[e^{e^k},se^{e^k})})=0.
\end{equation}
\begin{equation}\label{eqqq2}
\omega(\sum_k\chi_{(e^{k+e^k}/s,e^{k+e^k}]})=0.
\end{equation}
\end{lem}
\begin{proof} Denote the left hand side of \eqref{eqqq1} by $f(s).$ Due to the dilation invariance of $\omega,$ we have
$$f(s)=\omega(\sum_k\chi_{[te^{e^k},ste^{e^k})})=f(st)-f(t),\quad s,t>1.$$
Since $f$ is monotone and bounded, we have $f=0.$

Denote the left hand side of \eqref{eqqq2} by $g(s).$ Due to the dilation invariance of $\omega,$ we have
$$g(s)=\omega(\sum_k\chi_{(e^{k+e^k}/st,e^{k+e^k}/t]})=g(st)-g(t),\quad s,t>1.$$
Since $g$ is monotone and bounded, we have $g=0.$
\end{proof}

\begin{lem}\label{est12} Let $\omega$ be a dilation invariant generalised limit on $L_{\infty}(0,\infty).$ We have
\begin{enumerate}[i)]
\item $$\omega(\sum_k\frac{t}{\log(t)}e^{-e^k}\chi_{[e^{k-1+e^{k-1}},e^{k+e^k}]}(t))=0.$$
\item $$\omega(\sum_k\frac1{t\log(t)}e^{k+e^k}\chi_{[e^{e^k},e^{e^{k+1}}]}(t))=0.$$
\end{enumerate}
\end{lem}
\begin{proof} We only prove the first assertion. Proof of the second one is similar.

Fix $s>1.$ We have
$$\frac{t}{\log(t)}e^{-e^k}\leq\frac2s+2e^{-e^k/2},\quad \forall t\leq e^{k+e^k}/s,\quad \forall k\ge 1$$
and, therefore,
$$\sum_k\frac{t}{\log(t)}e^{-e^k}\chi_{[e^{k-1+e^{k-1}},e^{k+e^k}]}(t)\leq\frac2s+\sum_k\chi_{[e^{k+e^k}/s,e^{k+e^k}]}(t)+$$
$$+2\sum_ke^{-e^k/2}\chi_{[e^{k-1+e^{k-1}},e^{k+e^k}]}(t).$$
Clearly,
$$\omega(\sum_ke^{-e^k/2}\chi_{[e^{k-1+e^{k-1}},e^{k+e^k}]}(t))=0.$$
It follows from the Lemma \ref{est1} that
$$\omega(\sum_k\frac{t}{\log(t)}e^{-e^k}\chi_{[e^{k-1+e^{k-1}},e^{k+e^k}]}(t))\leq\frac2s.$$
Since $s$ is arbitrarily large, we have
$$\omega(\sum_k\frac{t}{\log(t)}e^{-e^k}\chi_{[e^{k-1+e^{k-1}},e^{k+e^k}]}(t))=0.$$
\end{proof}

\begin{lem}\label{omega construct} There exists a dilation invariant generalised limit $\omega$ on $L_{\infty}(0,\infty)$ such that
$$\omega(\sum_k\chi_{[e^{e^k},e^{k+e^k})})=1,\quad \omega(\sum_k\chi_{[e^{k+e^k},e^{e^{k+1}})})=0.$$
\end{lem}
\begin{proof} Define a positive, homogeneous functional $\pi$ on $L_{\infty}(0,\infty)$ by the formula
$$\pi(x)=\limsup_{N\to\infty}\frac1{\log(\log(N))}\int_N^{N\log(N)}x(s)\frac{ds}{s}.$$
It is verified in \cite[Lemma 4]{SSZ} that every $\omega\in L_{\infty}(0,\infty)^*$ satisfying $\omega\leq\pi$ is dilation invariant. Observing that
$$\pi(\sum_k\chi_{[e^{e^k},e^{k+e^k})})=1,$$
let us select $\omega\in L_{\infty}(0,\infty)^*$ satisfying $\omega\leq\pi$ and such that
$$\omega(\sum_k\chi_{[e^{e^k},e^{k+e^k})})=1.$$
Therefore,
$$\omega(\sum_k\chi_{[e^{k+e^k},e^{e^{k+1}})})=1-\omega(\sum_k\chi_{[e^{e^k},e^{k+e^k})})=0.$$
\end{proof}

Define a function $x$ by the formula
\begin{equation}\label{x counter}
x=\sup_{k\in\mathbb{N}}e^{-e^k}\chi_{[0,e^{k+e^k}]}.
\end{equation}
Fix $k\ge 1.$  For every $t\in[e^{k-1+e^{k-1}},e^{k+e^k}],$  we have
$$\frac1{\log(1+t)}\int_0^tx(s)ds\leq e^{1-k}\int_0^{e^{k+e^k}}x(s)ds\leq e^{1-k}\sum_{n=1}^ke^{-e^n}\cdot e^{n+e^n}\leq\frac{e^2}{e-1},$$
which guarantees $x\in M_{1,\infty}.$

\begin{lem}\label{est2} Let $x$ be as in \eqref{x counter} and let $\omega$ be as in Lemma \ref{omega construct}. We have $\tau_{\omega}(x)=(e-1)^{-1}.$
\end{lem}
\begin{proof} Fix $t\in[e^{k-1+e^{k-1}},e^{k+e^k}].$ We have
$$\int_0^tx(u)du=\frac{e^k}{e-1}+te^{-e^k}+O(1).$$
It follows that
$$\tau_{\omega}(x)=(e-1)^{-1}\omega(\sum_k\frac{e^k}{\log(t)}\chi_{[e^{k-1+e^{k-1}},e^{k+e^k}]}(t))+$$
$$+\omega(\sum_k\frac{t}{\log(t)}e^{-e^k}\chi_{[e^{k-1+e^{k-1}},e^{k+e^k}]}(t)).$$

By Lemma \ref{est12}, the second generalised limit above vanishes. We claim that the first generalised limit above is $1.$ Indeed,
$$\sum_k\frac{e^k}{\log(t)}\chi_{[e^{k-1+e^{k-1}},e^{k+e^k}]}(t)\geq(1+o(1))\sum_k\chi_{[e^{e^k},e^{k+e^k}]}(t)$$
and
$$\sum_k\frac{e^k}{\log(t)}\chi_{[e^{k-1+e^{k-1}},e^{k+e^k}]}(t)\leq\sum_k\chi_{[e^{e^k},e^{k+e^k}]}(t)+e\sum_k\chi_{[e^{k-1+e^{k-1}},e^{e^k}]}.$$
The claim follows from Lemma \ref{omega construct}.
\end{proof}

\begin{lem}\label{est3} Let $x$ be as in \eqref{x counter} and let $\omega$ be as in Lemma \ref{omega construct}. We have
$$\xi_{\omega}(x)=\frac{e}{e-1}\Gamma(1+\frac1q).$$
\end{lem}
\begin{proof} Fix $t\in[e^{e^k},e^{e^{k+1}}).$ We have
$$\int_{x>1/t}(x(u)-\frac1t)du=\frac{e^{k+1}}{e-1}-\frac1te^{k+e^k}+O(1).$$
This estimate and Lemma \ref{xiomega repr} yield
$$\frac1{\Gamma(1+1/q)}\xi_{\omega}(x)=\frac{e}{e-1}\omega(\sum_k\frac{e^k}{\log(t)}\chi_{[e^{e^k},e^{e^{k+1}}]}(t))-$$
$$-\omega(\sum_k\frac1{t\log(t)}e^{k+e^k}\chi_{[e^{e^k},e^{e^{k+1}}]}(t)).$$

It follows from Lemma \ref{est12} that the second generalised limit is $0.$ We claim that the first generalised limit is $1.$ Indeed,
$$\sum_k\frac{e^k}{\log(t)}\chi_{[e^{e^k},e^{e^{k+1}}]}(t)\geq(1+o(1))\sum_k\chi_{[e^{e^k},e^{k+e^k}]}$$
and
$$\sum_k\frac{e^k}{\log(t)}\chi_{[e^{e^k},e^{e^{k+1}}]}(t)\leq 1.$$
The claim follows from Lemma \ref{omega construct}.
\end{proof}

The following theorem delivers the promised counterexample.

\begin{thm}\label{counterexample} There exists $A\in\mathcal{M}_{1,\infty}$ and dilation invariant generalised limit $\omega$ on $L_{\infty}(0,\infty)$ such that
$$\Gamma(1+\frac1q)\tau_{\omega}(A)<\xi_{\omega}(A).$$
\end{thm}
\begin{proof} For brevity, we assume that the von Neumann algebra $\mathcal{N}$ is of type $II$ (the argument can be easily adjusted when $\mathcal{N}$ is of type $I$).  Let $x$ be as in \eqref{x counter} and let $A\in\mathcal{M}^+_{1,\infty}$ be such that $x=\mu(A).$ The assertion follows from Lemmas \ref{est2} and \ref{est3}.
\end{proof}

\section{Correctness of the definition for generalised heat kernel formulae}\label{bf def}

Let $\omega$ be a dilation invariant generalised limit on $L_{\infty}(0,\infty)$ and let $B\in\mathcal{N}.$ Following \cite{BF}, we consider the functionals on
$\mathcal{M}^+_{1,\infty}$ defined by the formula
\begin{equation}\label{xigeneral}
\xi_{\omega,B,f}(A)=(\omega\circ M)(t\to \frac1t\tau(f(tA)B)).
\end{equation}
The main result of this section, Theorem \ref{mhk bounded}, shows that the function
$$M\left(t\to \frac1t\tau(f(tA)B)\right)$$
is bounded, and so the formula \eqref{xigeneral} is well-defined.

\begin{lem}\label{operator estimate1} Let $A\in\mathcal{M}^+_{1,\infty}.$ We have $\tau(A^2e^A[0,1/t])=O(t^{-1}\log(t))$ as $t\to\infty.$
\end{lem}
\begin{proof} Let $c:=\|A\|_{1,\infty}.$ We have $\mu(s,A)\prec\prec c(1+s)^{-1}.$ Fix $t>0.$ Define decreasing function $x_t\in M_{1,\infty}(0,\infty)$ by setting
$$x_t(s)=
\begin{cases}
\frac{\log(1+ct\log(t))}{t\log(t)},& 0\leq s\leq ct\log(t)\\
\frac{c}{1+s}, & s>ct\log(t).
\end{cases}
$$
Define a decreasing function $y_t\in M_{1,\infty}(0,\infty)$ by setting
$$y_t(s)=\mu(A)\chi_{\{\mu(A)\leq 1/t\}}(s)+\frac1t\chi_{\{\mu(A)\geq1/t\}}(s),\ s>0.$$

We claim that $y_t\prec\prec x_t.$ Indeed, $y_t(s)\leq1/t\leq x_t(s)$ for $s\leq ct\log(t)$ and
$$\int_0^sy_t(u)du\leq c\int_0^s\frac{du}{1+u}=\int_0^sx_t(u)du$$
for $s>ct\log(t).$

It follows that
$$\tau(A^2e^A[0,\frac1t])\leq\int_0^{\infty}y_t^2(s)ds\leq\int_0^{\infty}x_t^2(s)ds.$$
We have
$$\int_0^{\infty}x_t^2(s)ds=\frac{c\log^2(1+ct\log(t))}{t\log(t)}+\int_{ct\log(t)}^{\infty}\frac{c^2}{(1+s)^2}ds\leq 5c\frac{\log(t)}{t}.$$
\end{proof}

\begin{lem}\label{square estimate} Let $f(t)=t^2\chi_{[0,1]}(t)$ and let $A\in\mathcal{M}^+_{1,\infty}.$ We have
$$t\to M(\frac1t\tau(f(tA)))\in L_{\infty}(0,\infty).$$
\end{lem}
\begin{proof} For fixed $t>0,$ we have
$$M(\frac1t\tau(f(tA)))=\frac1{\log(t)}\int_1^t\tau(A^2e^A[0,\frac1s])ds=\frac1{\log(t)}\tau(A^2\int_1^te^A[0,\frac1s]ds).$$

Integrating by parts, we obtain
$$\int_1^te^A[0,\frac1s]ds=se^A[0,\frac1s]|_1^t-\int_1^tsde^A[0,\frac1s]=se^A[0,\frac1s]|_1^t+\int_{1/t}^1u^{-1}de^A[\frac1t,u]=$$
$$=O(1)+A^{-1}e^A[\frac1t,\infty]+te^A[0,\frac1t].$$
Therefore,
$$M(\frac1t\tau(f(tA)))=\frac1{\log(t)}\tau(Ae^A(\frac1t,\infty))+\frac{t}{\log(t)}\tau(A^2e^A[0,\frac1t])+O(\frac1{\log(t)}).$$
It follows from the definitions of $\|\cdot\|_{1,\infty}$ and $d_A(\cdot)$ that for every $A\in\mathcal{M}_{1,\infty}$ and every $t>0,$ we have
$$d_A(\frac1t)\leq\max\{1,\|A\|_{1,\infty}\}\log(1+t).$$
Clearly,
$$\frac1{\log(t)}\tau(Ae^A[0,\frac1t])=\frac1{\log(t)}\int_0^{d_A(1/t)}\mu(s,A)ds\leq\frac{\log(d_A(1/t))}{\log(t)}\|A\|_{1,\infty}\in L_{\infty}.$$
The assertion follows now from the Lemma \ref{operator estimate1}.
\end{proof}

\begin{thm}\label{mhk bounded} Let a bounded function $f\in C^2[0,\infty)$ be such that $f(0)=f'(0)=0$. Let $A\in\mathcal{M}^+_{1,\infty}$ and let $B\in\mathcal{N}.$ We have
$$M\left(t\to \frac1t\tau(f(tA)B)\right)\in L_{\infty}(0,\infty).$$
\end{thm}
\begin{proof} Due to the well known inequality $\tau(CB)\leq\tau(|C|)\|B\|,$ it suffices to prove the theorem only when $B=1.$ In this case, for the function $f(t):=t^2\chi_{[0,1]}(t),$ the assertion follows from Lemma \ref{square estimate}. If $f(t):=\chi_{(1,\infty)}(t)$ then it holds trivially. Thus, it holds  for the function $f(t):=\min\{1,t^2\}.$ Finally, observe that the assumptions on $f$ guarantee that there exists a constant $c>0$ such that $|f(t)|\leq c\min\{1,t^2\}.$
\end{proof}

Since the function $t\to\exp(-t^{-q})$ satisfies the assumptions of Theorem \ref{mhk bounded} we obtain the following corollary, which was implicitly proved in \cite[Theorem 5.2]{CRSS}.

\begin{cor} For every $q>0$ and every $A\in\mathcal{M}^+_{1,\infty},$ we have
$$M\left(t\to \frac1t\tau(\exp(-(tA)^{-q}))\right)\in L_{\infty}(0,\infty).$$
\end{cor}

\section{Reduction theorem for generalised heat kernel formulae}\label{bf red}

The results of this section extend and generalise those of \cite[Theorem 4.1]{CPS} and \cite[Theorem 5.2]{CRSS}. We also give an answer to the question asked in \cite[page 52]{BF}. We explicitly prove that the functional $\xi_{\omega,B,f}$ (extended to  $\mathcal{M}_{1,\infty}$ as in \eqref{xi-full}) is linear on $\mathcal{M}_{1,\infty}.$

\begin{lem}\label{golova k nulyu} Let $f\in C^2[0,\infty)$ be such that $f(0)=f'(0)=0.$ Let $A\in\mathcal{M}^+_{1,\infty}$ and let $B\in \mathcal{N}.$ For every dilation invariant generalised limit $\omega$ on $L_\infty(0,\infty),$ we have
$$\lim_{\varepsilon\to0}(\omega\circ M)(\frac1t\tau(f(tAe^A[0,\frac{\varepsilon}{t}])B))=0.$$
\end{lem}
\begin{proof} Since $|f(t)|\leq const\cdot t^2$ for $t\in[0,1],$ it is sufficient to prove the assertion for $f(t)=t^2.$ As in the proof of Theorem \ref{mhk bounded}, it is sufficient to assume that $B=1.$

By Theorem \ref{mhk bounded}, for every $\varepsilon>0$ we have
$$M\left(t\to \frac1t\tau((tAe^A[0,\frac{\varepsilon}{t}])^2)\right)\in L_{\infty}(0,\infty).$$
Since $\omega$ is dilation invariant, we conclude
$$(\omega\circ M)(\frac1t\tau((tAe^A[0,\frac{\varepsilon}{t}])^2))=\varepsilon(\omega\circ M)(\frac1t\tau((tAe^A[0,\frac1t])^2)).$$
The assertion follows immediately.
\end{proof}

\begin{lem}\label{hvost k nulyu} Let $f\in L_{\infty}(0,\infty)$ be such that $f(0)=0.$ Let $A\in\mathcal{M}^+_{1,\infty}$ and let $B\in\mathcal{N}.$ For every dilation invariant generalised limit $\omega$ on $L_\infty(0,\infty),$ we have
$$\lim_{\varepsilon\to 0}(\omega\circ M)(\frac1t\tau(f(tAe^A(\frac{1}{\varepsilon t},\infty))B))=0.$$
\end{lem}
\begin{proof} As before, we may assume that $B=1.$ It is clear that
$$f(tAe^A(\frac{1}{\varepsilon t},\infty))\leq\|f\|e^A(\frac{1}{\varepsilon t},\infty).$$

Since $\omega\circ M$ is dilation invariant, we obtain
$$(\omega\circ M)(\frac1t\tau(e^A(\frac{1}{\varepsilon t},\infty)))=\varepsilon(\omega\circ M)(\frac1td_A(\frac1t)).$$
The assertion follows immediately.
\end{proof}

\begin{lem}\label{seredina} Let $f:\mathbb{R}_+\to\mathbb{R}$ be monotone on $[a,b]$ and such that $f(0)=0.$ Let $A\in\mathcal{M}^+_{1,\infty}$ and let $B\in\mathcal{N}.$ For every dilation invariant generalised limit $\omega$ on $L_\infty(0,\infty)$ we have
$$(\omega\circ M)(\frac1t\tau(f(tAe^A[\frac{a}{t},\frac{b}{t}))B))=(\int_a^bf(s)\frac{ds}{s^2})\cdot(\omega\circ M)(\frac1t\tau(e^A[\frac1t,\infty)B)).$$
\end{lem}
\begin{proof} Without loss of generality, we may assume that $f$ is increasing on $[a,b]$ and that $B\ge 0.$

Let $a=a_0\leq a_1\leq a_2\leq\cdots\leq a_n=b.$ For every given $t>0,$ we have
$$e^A[\frac{a}{t},\frac{b}{t})=\sum_{k=0}^{n-1}e^A[\frac{a_k}t,\frac{a_{k+1}}t).$$

Since $f$ is increasing on $[a,b]$ and $f(0)=0,$ we have
$$f(a_k)e^A[\frac{a_k}t,\frac{a_{k+1}}t)\leq f(tAe^A[\frac{a_k}t,\frac{a_{k+1}}t))\leq f(a_{k+1})e^A[\frac{a_k}t,\frac{a_{k+1}}t).$$
Therefore,
$$(\omega\circ M)(\frac1t\tau(f(tAe^A[\frac{a}{t},\frac{b}{t}))B))\leq\sum_{k=0}^{n-1}f(a_{k+1})(\omega\circ M)(\frac1t\tau(e^A[\frac{a_k}t,\frac{a_{k+1}}t)B))$$
and
$$(\omega\circ M)(\frac1t\tau(f(tAe^A[\frac{a}{t},\frac{b}{t}))B))\geq\sum_{k=0}^{n-1}f(a_k)(\omega\circ M)(\frac1t\tau(e^A[\frac{a_k}t,\frac{a_{k+1}}t)B)).$$
We have
$$e^A[\frac{a_k}t,\frac{a_{k+1}}t)=e^A[\frac{a_k}t,\infty)-e^A[\frac{a_{k+1}}t,\infty).$$
For all $c>0,$ we have
$$(\omega\circ M)(\frac1t\tau(e^A(\frac{c}{t},\infty)B))=c^{-1}(\omega\circ M)(\frac1t\tau(e^A(\frac1t,\infty)B)).$$
Therefore,
$$(\omega\circ M)(\frac1t\tau(e^A[\frac{a_k}t,\frac{a_{k+1}}t)B))=(\frac1{a_k}-\frac1{a_{k+1}})(\omega\circ M)(\frac1t\tau(e^A(\frac1t,\infty)B)).$$
Hence,
$$(\sum_{k=0}^{n-1}f(a_k)(\frac1{a_k}-\frac1{a_{k+1}}))(\omega\circ M)(\frac1t\tau(e^A(\frac1t,\infty)B))\leq$$
$$\leq(\omega\circ M)(\frac1t\tau(f(tAe^A[\frac{a}{t},\frac{b}{t}))B))\leq$$
$$\leq(\sum_{k=0}^{n-1}f(a_{k+1})(\frac1{a_k}-\frac1{a_{k+1}}))(\omega\circ M)(\frac1t\tau(e^A(\frac1t,\infty)B)).$$
Both coefficients in the latter formula tend to $\int_a^bf(s)s^{-2}ds.$
\end{proof}

\begin{lem}\label{sborka} Let a bounded function $f\in C^2[0,\infty)$ be such that $f(0)=f'(0)=0.$ Let $A\in\mathcal{M}^+_{1,\infty}$ and let $B\in\mathcal{N}.$ For every dilation invariant generalised limit $\omega$ on $L_\infty(0,\infty)$  we have
$$\xi_{\omega,B,f}(A)=(\int_0^{\infty}f(s)\frac{ds}{s^2})(\omega\circ M)(\frac1t\tau(e^A(\frac1t,\infty)B)).$$
\end{lem}
\begin{proof} Let $f$ satisfy the assumptions above. Observe that the assertion of Lemma \ref{seredina} holds for the function $f|_{[a,b]},$ where $0<a<b<\infty.$ Indeed, every such function is a function of bounded variation and therefore may be written as a difference of two monotone functions. Now the assertion follows from Lemmas \ref{golova k nulyu},\ref{hvost k nulyu},\ref{seredina} by setting $a:=\varepsilon$ and $b:=\varepsilon^{-1}$ and letting $\varepsilon\to 0.$
\end{proof}

\begin{cor}\label{corollary 46} Let a bounded function $f\in C^2[0,\infty)$ be such that $f(0)=f'(0)=0.$ Let $A\in\mathcal{M}^+_{1,\infty}$ and let $B\in\mathcal{N}^+.$ For every dilation invariant generalised limit $\omega$ on $L_\infty(0,\infty)$  we have
$$\xi_{\omega,B,f}(A)=(\int_0^{\infty}f(s)\frac{ds}{s^2})\omega(\frac1{\log(1+t)}\tau((A-\frac1t)e^A(\frac1t,\infty)B)).$$
\end{cor}
\begin{proof} It follows from the definition of Cesaro operator $M$ that
$$M\left(t\to\frac1t\tau(e^A(\frac1t,\infty)B)\right)=\frac1{\log(t)}\int_1^t\tau(e^A(\frac1s,\infty)B)\frac{ds}{s^2}.$$
Integrating by parts, we obtain
$$\frac1{\log(t)}\int_1^t\tau(e^A(\frac1s,\infty)B)\frac{ds}{s^2}=\frac1{\log(t)}\int_{1/t}^1\tau(e^A(u,\infty)B)du=$$
$$=\frac1{\log(t)}\cdot u\tau(e^A(u,\infty)B)|_{1/t}^1-\frac1{\log(t)}\int_{1/t}^1ud\tau(e^A(u,\infty)B)=$$
$$=\frac{-1}{t\log(t)}\cdot\tau(e^A(\frac1t,\infty)B)+\frac{-1}{\log(t)}\tau(\int_{1/t}^{\infty}ude^A(u,\infty)B)+o(1).$$
Evidently,
$$-\tau(\int_{1/t}^{\infty}ude^A(u,\infty)B)=\tau(Ae^A(\frac1t,\infty)B).$$
Therefore,
$$M\left(t\to\frac1t\tau(e^A(\frac1t,\infty)B)\right)=\frac1{\log(t)}\tau((A-\frac1t)e^A(\frac1t,\infty)B)+o(1).$$
The assertion follows now from Lemma \ref{sborka}.
\end{proof}

The first assertion in lemma below can be found in \cite[Theorem 11]{Brown-Kosaki}. For the second assertion we refer to \cite[Theorem 3.5]{Bourin}.

\begin{lem}\label{main heat estimate} Let $A,B\in B^+(H)$ and let $f$ be convex continuous function such that $f(0)=0.$ We have
\begin{enumerate}[i)]
\item $\tau(B^{1/2}f(A)B^{1/2})\geq\tau(f(B^{1/2}AB^{1/2}))$ if $B\leq 1.$
\item $\tau(B^{1/2}f(A)B^{1/2})\leq\tau(f(B^{1/2}AB^{1/2}))$ if $B\geq 1.$
\end{enumerate}
\end{lem}

We show in the following lemma that $\xi_{\omega,B,f}$ depends continuously on $B.$

\begin{lem}\label{hk continuity} If $A\in\mathcal{M}^+_{1,\infty}$ and let $B_n, B\in \mathcal{N},$ $n\ge 1,$ then
$$\|\xi_{\omega,B_n}(A)-\xi_{\omega,B}(A)\|\leq\xi_{\omega}(A)\cdot\|B_n-B\|.$$
\end{lem}
\begin{proof} The assertion follows from the inequality
$$|\tau(f(tA)B_n)-\tau(f(tA)B)|\leq\tau(f(tA))\cdot\|B_n-B\|.$$
\end{proof}

The following theorem extends the results of \cite{CPS,CRSS} and gives an affirmative answer to the question stated in \cite{BF}. It also shows that the functionals $\xi_{\omega,B,f}(\cdot)$ are linear functionals on $\mathcal{M}_{1,\infty}$ for a wide class of functions $f.$

\begin{thm}\label{xihk red} Let a bounded function $f\in C^2[0,\infty)$ be such that $f(0)=f'(0)=0.$ Let $A\in\mathcal{M}_{1,\infty}$ and let $B\in\mathcal{N}.$ For every dilation invariant generalised limit $\omega$ on $L_\infty(0,\infty)$   we have
\begin{equation}\label{BF-heat reduction formula}
\xi_{\omega,B,f}(A)=\frac1{\Gamma(1+1/q)}(\int_0^{\infty}f(s)\frac{ds}{s^2})\xi_{\omega}(AB).
\end{equation}
\end{thm}
\begin{proof} It follows from Theorem \ref{xiomega linearity} that $\xi_{\omega}$ is linear and fully symmetric. By Theorem \ref{fs is dixmier} and \eqref{dixmier trace pr2}), we have $\xi_{\omega}(B^{1/2}AB^{1/2})=\xi_{\omega}(AB).$

Recall that function $u\to(u-1/t)_+$ is convex. It follows from Lemma \ref{main heat estimate} that
\begin{enumerate}[i)]
\item $\tau((A-\frac1t)_+B)\geq\tau((B^{1/2}AB^{1/2}-\frac1t)_+)$ if $B\leq 1.$
\item $\tau((A-\frac1t)_+B)\leq\tau((B^{1/2}AB^{1/2}-\frac1t)_+)$ if $B\geq 1.$
\end{enumerate}

It follows from Corollary \ref{corollary 46} that for $0\leq B\leq1$ we have
\begin{equation}\label{hk lower estimate}
\xi_{\omega,B,f}(A)\geq\frac1{\Gamma(1+1/q)}(\int_0^{\infty}f(s)\frac{ds}{s^2})\xi_{\omega}(B^{1/2}AB^{1/2}).
\end{equation}
Since both sides are homogeneous, the inequality \eqref{hk lower estimate} is valid for every $B.$

It follows from \ref{corollary 46} that for $B\geq1$ we have
\begin{equation}\label{hk upper estimate}
\xi_{\omega,B,f}(A)\leq\frac1{\Gamma(1+1/q)}(\int_0^{\infty}f(s)\frac{ds}{s^2})\xi_{\omega}(B^{1/2}AB^{1/2}).
\end{equation}
Since both sides are homogeneous, the inequality \eqref{hk upper estimate} is valid if $B$ is bounded from below by a strictly positive constant.

Thus, we have the equality \eqref{BF-heat reduction formula} valid for every $B$ bounded from below by a strictly positive constant. Set $B_n=Be^B(1/n,\infty)+1/ne^B[0,1/n].$ It follows that equality \eqref{BF-heat reduction formula} holds with $B$ replaced with $B_n$ throughout. By Lemma \ref{hk continuity}, we have $\xi_{\omega,B_n,f}(A)\to\xi_{\omega,B,f}(A).$ Since $AB_n\to AB$ in $\mathcal{M}_{1,\infty}$ and since $\xi_{\omega}$ is bounded on $\mathcal{M}_{1,\infty},$ we have $\xi_{\omega}(AB_n)\to\xi_{\omega}(AB).$ The assertion follows immediately.
\end{proof}

The following corollary treats the case of classical heat kernel formulae. We use the notation
$$\xi_{\omega,B}(A)=(\omega\circ M)(\frac1t\tau(\exp(-(tA)^{-q})B)).$$

\begin{cor} Let $A\in\mathcal{M}^+_{1,\infty}$ and let $B\in\mathcal{N}.$ For every dilation invariant generalised limit $\omega$ on $L_\infty(0,\infty)$   we have $\xi_{\omega,B}(A)=\xi_{\omega}(AB).$
\end{cor}
\begin{proof} Use $f(t)=\exp(-t^{-q})$ in Theorem \ref{xihk red} and observe that
$$\int_0^{\infty}f(s)\frac{ds}{s^2}=\Gamma(1+\frac1q).$$
\end{proof}

The following assertion extends \cite[Theorem 33]{SSZ}.

\begin{cor}\label{corollary 51} Let $A\in\mathcal{M}^+_{1,\infty}$ and let $B\in\mathcal{N}.$ For every dilation invariant generalised limit $\omega$ on $L_\infty(0,\infty)$ such that $\omega=\omega\circ M,$ we have
$$\xi_{\omega,B}(A)=\Gamma(1+\frac1q)\tau_{\omega}(AB).$$
\end{cor}

\end{document}